\documentclass[12pt]{article}
 
\usepackage{amssymb,amsfonts,amsmath,latexsym,epsf,tikz,url}
\usepackage[usenames,dvipsnames]{pstricks}
\usepackage{pstricks-add}
\usepackage{epsfig}
\usepackage{pst-grad} 
\usepackage{pst-plot} 
\usepackage[space]{grffile} 
\usepackage{etoolbox} 
\usepackage{float}
\usepackage{soul}
\usepackage{tikz}
\usepackage[colorinlistoftodos]{todonotes}
\usepackage{pgfplots}
\usepackage{mathrsfs}
\usepackage[colorlinks]{hyperref}
\usetikzlibrary{arrows}
\makeatletter 
\patchcmd\Gread@eps{\@inputcheck#1 }{\@inputcheck"#1"\relax}{}{}
\makeatother

\newtheorem{theorem}{Theorem}[section]

\newtheorem{proposition}[theorem]{Proposition}
\newtheorem{observation}[theorem]{Observation}
\newtheorem{corollary}[theorem]{Corollary}
\newtheorem{lemma}[theorem]{Lemma}
\newtheorem{remark}[theorem]{Remark}

\newtheorem{definition}[theorem]{Definition}

\tikzstyle{black_v}=[fill=black, draw=black, shape=circle]
\tikzstyle{none}=[fill=none, draw=none, shape=circle]
\tikzstyle{blue_v}=[fill=blue, draw=blue, shape=circle]
\tikzstyle{red_v}=[fill=red, draw={rgb,255: red,246; green,10; blue,34}, shape=circle]
\tikzstyle{green_v}=[fill={rgb,255: red,17; green,255; blue,0}, 
	draw={rgb,255: red,45; green,255; blue,8}, shape=circle]
\tikzstyle{BigBlue}=[fill=blue, draw=blue, shape=circle, scale=1.3]
\tikzstyle{BigRed}=[fill=red, draw=red, shape=circle, scale=1.75]
\tikzstyle{BBigBlue}=[fill=blue, draw=blue, shape=circle, scale=1.75]
\tikzstyle{BigGreen}=[fill={rgb,255: red,49; green,215; blue,37}, 
	draw={rgb,255: red,0; green,184; blue,0}, shape=circle]

\tikzstyle{red_E}=[-, draw=red, fill=red, ultra thick]
\tikzstyle{dashed_line}=[-, dashed]
\tikzstyle{green_E}=[-, draw={rgb,255: red,58; green,228; blue,83}]
\tikzstyle{magenta_E}=[-, draw={rgb,255: red,246; green,101; blue,246}]
\tikzstyle{blue_E}=[-, draw={rgb,255: red,32; green,32; blue,253}, ultra thick]
\tikzstyle{olive_E}=[-, draw={rgb,255: red,0; green,128; blue,128}]
\tikzstyle{flecha}=[->]
\tikzstyle{doble}=[-, double]
\tikzstyle{dots}=[-, dotted, tikzit draw={rgb,255: red,238; green,87; blue,236}]
\tikzstyle{gray_e}=[-, fill=none, draw={rgb,255: red,171; green,171; blue,171}]
\tikzstyle{blue_e}=[-, draw={rgb,255: red,28; green,93; blue,244}]

\begin{document}

\def\nt{\noindent}

\title{Total restrained coalitions in graphs}

\bigskip 
\author{
 M. Chellali $^{1}$,  
J.C. Valenzuela-Tripodoro $^{2}$,  
H. Golmohammadi $^{3,4}$,  \\[.5em]
I.I. Takhonov $^{3}$,
N.A. Matrokhin $^{3}$  
}


\maketitle
\begin{center}
    
$^{1}$LAMDA-RO Laboratory, Department of Mathematics, University 
of Blida, Blida, Algeria

$^{2}$Department of Mathematics, University of C\'{a}diz, Spain

$^{3}$Novosibirsk State University, Pirogova str. 2, Novosibirsk, 
630090, Russia\\ 

\medskip
$^{4}$Sobolev Institute of Mathematics, Ak. Koptyug av. 4, 
Novosibirsk,
630090, Russia\\

\medskip
    	{\tt m\_chellali@yahoo.com ~~ jcarlos.valenzuela@uca.es ~~
    	h.golmohammadi@g.nsu.ru ~~ i.takhonov@g.nsu.ru ~~ 
    	n.matrokhin@g.nsu.ru} 

\end{center}

\begin{abstract}

In this paper, we introduce the concept of total restrained coalition 
and total restrained coalition partition in graphs. A vertex set in a graph 
without isolated vertices is a total restrained dominating 
set (TRD-set) if it is dominating, induces a subgraph without isolated 
vertices, and the vertices not in the set also induce a subgraph 
without isolated vertices. Two vertex sets, which are not TRD-sets, 
form a total restrained coalition if their union is a TRD-set. A total 
restrained coalition partition is a partition where none of its elements 
are TRD-sets, but each forms a total restrained coalition with 
another element. The goal is to maximize the cardinality of such a partition, 
denoted $C_{tr}(G)$. We initiate the study of this concept by proving 
certain properties, extremal values, general bounds, 
and its relation to known structural parameters. Exact values for specific 
graph families are also provided.
\end{abstract}

\noindent{\bf Keywords:}   Coalition; total restrained coalition, 
total restrained dominating set.
  
\medskip
\noindent{\bf AMS Subj.\ Class.:}  05C60. 


\section{Introduction} 

Throughout this article, we only consider finite and simple graphs 
without isolated vertices. For such a graph $G=(V,E)$ and a vertex 
$v\in V$, we denote by $N(v):= \{w\in V\mid vw\in E\}$ the open 
neighborhood of $v$ and by $N[v] := N(v)\cup\{v\}$ its closed 
neighborhood.
The order of a graph $G$ refers to the cardinality $|V|$ of its 
set of vertices. Each vertex of $N(v)$ is called a neighbor
of $v$, and the cardinality of $|N(v)|$ is called the degree of $v$, 
denoted by $deg(v)$. 
The minimum and maximum degree of graph vertices are denoted by 
$\delta(G)$ and $\Delta(G)$, respectively. An isolated vertex in $G$ 
is a vertex of degree 0. A graph is isolate-free if it contains no 
isolated vertex.

A set $S \subseteq V$ is called a dominating set if every vertex of 
$V \setminus S$ is adjacent to at least one vertex in $S$. Further 
if every vertex in $G$ is adjacent to some other vertex in $S$, 
then $S$ is a total dominating set, abbreviated TD-set of $G$. The 
domination number of $G$, denoted by $\gamma(G)$, is the minimum 
cardinality of a dominating set of $G$, while the total domination 
number $\gamma_{t}(G)$ of $G$ is the minimum cardinality of a TD-set 
of $G$. Various aspects of domination are well studied in the literature, 
and a thorough study of domination appears in \cite{A11, A12}.
  
Given a graph $G$, a set $S \subseteq V (G)$ is said to be a total 
restrained dominating set (abbreviated TRD-set) of $G$ if every vertex 
in $V\setminus S$ is adjacent to at least one vertex in $S$ and at least 
one other vertex in $V\setminus S$, and every vertex in $S$ is adjacent 
to at least one other vertex in $S$.

The total restrained domination number of $G$,
denoted by $\gamma_{tr}(G)$, is the cardinality of a minimum TRD-set of 
$G$. It is worth mentioning that every graph without isolated vertices 
has a TRD-set, since $S=V$ is such a set. The concept of the total 
restrained domination was introduced by Telle and Proskurovsky \cite{A14}, 
although implicitly, as a vertex partitioning problem. Total restrained 
domination in graphs is well studied in the literature. For more details 
we refer the reader to the recent book chapter by Hattingh and 
Joubert \cite{A6}.

Let $\mathcal{D}$ be a partition of the vertex set $V(G)$ of $G$. If 
all sets of $\mathcal{D}$ are total dominating sets in $G$, then 
$\mathcal{D}$ is called a total domatic partition of $G$. The maximum 
number of sets of a total domatic partition of $G $ is the total domatic 
number $d_{t}(G)$ of $G$. In \cite{Z}, Zelinka studied this concept. 
Analogously the total restrained domatic partition is a partition
of vertices of a graph into total restrained dominating sets. The maximum 
cardinality of a total restrained domatic partition is called the total 
restrained domatic number, denoted by $d^r_{t}(G)$. The total restrained 
domatic number of a graph was introduced by Zelinka in \cite{A15}.\newline 
Fairly recently, the concept of coalition in graphs has triggered a great 
deal of interest due to its definition, which is based on dominating sets.
A coalition in a graph $G$ is composed of two disjoint sets of vertices 
$X$ and $Y$ of $G$, neither of which is a dominating set but whose union 
$X \cup Y$ is a dominating set of $G $. A coalition partition is a vertex 
partition $\pi=\{V_1,V_2,\dots,V_k\}$ of $V$ such that for every 
$i\in\{1,2,\dots,k\}$ the set $V_i$ is either a dominating set and 
$|V_i|=1$, or there exists another set $V_j$ so that they form a coalition.
The maximum cardinality of a coalition partition is called the coalition 
number of the graph, and denoted by $C(G)$. Coalitions in graphs were 
introduced and first studied by Haynes et al. in \cite{A7}, and have been 
studied further \cite{A4,A8,A9,A10}. Several types of domination coalitions 
have been studied by imposing additional conditions on the domination 
coalition, see \cite{A1,A2,A3,A5,A13}. The aim of this paper is to introduce 
and study the concept of total restrained coalition in graphs. We begin with
the following definitions.

\begin{definition}[Total restrained coalition]
Two disjoint sets $X,Y\subseteq V(G)$ form a total restrained coalition 
in a graph $G$ if they are not TRD-sets but their union is a TRD-set in $G$. 
 \end{definition} 
 
\begin{definition}[Total restrained coalition partition]\label{2.2} 
A total restrained coalition partition, {  abbreviated as a 
trc-partition}, 
of a graph $G$ is a partition $\Phi=\{V_1,V_2,\dots,V_k\}$ of the vertex 
set $V$ such that any $V_i\in \Phi, 1\leq i \leq k,$ is not a TRD-set but 
forms a total restrained coalition with another {  set 
$V_j \in \Phi$}. The maximum cardinality of a total restrained coalition 
partition is called the total restrained coalition number of $G$
and denoted by $C_{tr}(G)$. A trc-partition of $G$ of cardinality $C_{tr}(G)$
is called a $C_{tr}(G)$-partition.\medskip
\end{definition}

Since every TRD-set in \(G\) is a TD-set, a natural question 
that arises is whether both problems are equivalent.
Consider the cycle graph \(C_3\) with \(V(C_3) = \{x, y, z\}\). The 
trc-partitions of \(C_3\) with two elements are  
\[
\Phi_1 = \{\{x\}, \{y, z\}\}, \quad \Phi_2 = \{\{y\}, \{x, z\}\}, \quad 
\Phi_3 = \{\{z\}, \{x, y\}\}.
\]  
First, note that none of the trc-partitions \(\Phi_1\), \(\Phi_2\),
 or \(\Phi_3\) qualifies as a tc-partition, because each contains a 
 two-vertex set that is a total dominating set.

Furthermore, it is straightforward to see that \(\{\{x\}, \{y\}, \{z\}\}\) 
is a tc-partition but not a trc-partition of \(C_3\), leading 
to the inequality  
\[
2 = C_{\text{tr}}(C_3) < C_t(C_3) = 3.
\]
Therefore, both problems are not equivalent, and it is worth studying 
the total restrained coalition partition problem.

\medskip 
The main contributions of this work are as follows. In Section 2, 
we first discuss the possibility of the existence of trc-partitions 
in graphs and derive some bounds. In Section 3, we determine the 
total restrained coalition number for some classes of graphs. In 
Section 4, we are interested in graphs with a large total restrained 
coalition number.

\section{Properties and bounds}

In this section, we present basic properties and bounds on the total 
restrained coalition number. We first call up the following trivial 
observation that we need for what follows.

\begin{observation}{\rm\cite{A6}}
    Every graph $G$ without an isolated vertex has a TRD-set.
\end{observation}

Now we state the following observation about the total restrained 
coalition number of a graph $G$.

\begin{observation} 
   If a graph $G$ contains an isolated vertex, then $C_{tr}(G)=0$.
\end{observation}
    
We are now in a position to prove the following result.

\begin{theorem} ~\label{1}
Let $G$ be an isolate-free graph. Then $G$ has, at least, a 
trc-partition and $C_{tr}(G)\ge 2d_t^r(G).$
\end{theorem}

\begin{proof} Consider a graph $G$ with a total restrained domatic 
partition $\mathcal{D}=\{S_1, \ldots, S_k\}$, {  with $k=d_t^r$}. 
In what follows we demonstrate the process of constructing a trc-partition 
$\Phi$ of $G$. For any integer $1\leq i\leq k-1$, assume that $S_i$ is a 
minimal TRD-set of $G$. If it is not, then there exists a minimal 
TRD-set $S_i'\subseteq S_i$. In this case, we replace $S_i$ with 
$S_i'$ and put all members of $S_i\setminus S_i'$ to $S_k$. In order 
to create a {  trc-partition, $\Phi$ of $G$,} we divide 
each minimal TRD-set $S_i$ with $i<k$ into two non-empty sets $S_{i,1}$ 
and $S_{i,2}$ and add them to $\Phi$. Note that neither $S_{i,1}$ 
nor $S_{i,2}$ is a TRD-set, but their union is a TRD-set. Next, we 
consider the set $S_k$. If $S_k$ is a minimal TRD-set, we split it 
into two non-empty sets $S_{k,1}$ and $S_{k,2}$ and attach them to 
$\Phi$. So, we obtain a trc-partition $\Phi$ of cardinality $2d_t^r.$

If $S_k$ is not a minimal TRD-set, there exists a set $S_k'\subseteq S_k$ 
that is minimal and total restrained dominating. We split $S_k'$ into 
two non-empty sets $S_{k,1}'$ and $S_{k,2}'$ and attach them to $\Phi$. 
Let $S_k''=S_k\backslash S_k'$. It is worth emphasizing that $S_k''$ 
cannot be a TRD-set, as this would imply that $d^r_{t}(G)>k$, 
against our assumptions. 
If $S_k''$ forms a total restrained coalition with any set in $\Phi$, 
we attach it to $\Phi$ and finish the construction process 
obtaining a total restrained coalition partition 
$\Phi$, of cardinality at least $2k+1\ge 2d_t^r$. Otherwise, by 
replacing $S_{k,2}'$ with $S_{k,2}'\cup S_k''$ in $\Phi$ we obtain 
a trc-partition with cardinality $2k=2d_t^r$. 
$\Box$

\end{proof}\medskip

It is clear that for all graphs $G$ without isolated vertices, $d^r_{t}(G)\geq 1$. 
By Theorem \ref{1} we infer the following result.

\begin{corollary}~\label{3}
 If $G$ is an isolate-free graph, then $2\leq C_{tr}(G)\leq n$.
\end{corollary}

Notice that if an isolate-free graph $G$ satisfies $C_{tr}(G)=2$, then we must have 
$d^r_{t}(G)=1$. However, the converse is not true and this can be seen by the cycle $C_5$, 
where $d^r_{t}(C_5)=1$ and $C_{tr}(C_5)=3$.\medskip

We next recall the following result due to Zelinka \cite{A15}.

\begin{theorem} {\rm\cite{A15}}\label{A}
Let $G$ be a graph without isolated vertices. Then $d^r_{t}(G)=d_{t}(G)$.
\end{theorem}

Plugging the result of Theorem \ref{A} into the bound of Theorem \ref{1} 
immediately yields the following result.

\begin{corollary}\label{B}
    Let $G$ be a graph without isolated vertices. Then $C_{tr}(G)\geq 2d_{t}(G)$.
\end{corollary}

In \cite{Z}, Zelinka showed that if $G$ is an isolate-free graph of order $n$ and 
minimum degree $\delta$,
then $d_{t}(G)\geq\left\lfloor \frac{n}{n-\delta+1}\right\rfloor$. As a
consequence of this result and Corollary \ref{B}, we have the following result.

\begin{corollary}
\label{delta}For any isolate-free graph $G,$ 
$C_{tr}(G)\geq2\left\lfloor \frac{n}{n-\delta+1}\right\rfloor$.
\end{corollary}

Restricted to connected graphs $G$ with minimum degree at least two and girth seven or more, 
we provide a lower bound for $C_{tr}(G)$ in terms of the maximum degree.

\begin{theorem} \label{girth 7} 
Let $G$ be a connected graph with minimum degree $\delta
(G)\geq2,$ maximum degree $\Delta(G)$ and girth at least $7.$ Then
$C_{tr}(G)\geq\Delta(G)+1.$
\end{theorem}

\textbf{Proof. }Let $\delta(G)=\delta$ and $\Delta(G)=\Delta.$ Let $w$ be a
vertex with maximum degree, and let $w_{1},w_{2},...,w_{\Delta}$ denote the
neighbors of $w$. Clearly, $N(w)$ is independent, for otherwise $G$ has a
triangle contradicting the assumption on the girth. The same argument of the
girth together with the fact $\delta\geq2$ also imply $V(G)-N[w]$ is non
empty. Let $A=V(G)-N(w).$ Clearly, since $\delta\geq2,$ each
$w_{i}\in N(w)$ has at least one neighbor in $A$ other than $w.$ 
For any $w_{i}\in N(w),$ let $w_{i}^{\prime}$ denote a neighbor of $w_{i}$ in
$A-\{w\}.$ Recall that $w$ has no neighbor in $A$ and
thus $ww_{i}^{\prime}\notin E(G).$ We make some useful remarks for the following. For any 
two distinct vertices $w_{i},w_{j}\in N(w),$ we have:

(i)
$w_{i}^{\prime}\neq w_{j}^{\prime}$, for otherwise vertices $w,w_{i},w_{j}$
and $w_{i}^{\prime}$ induce a cycle $C_{4},$ contradicting $G$ has girth at
least 7.

(ii)
$w_{i}^{\prime}w_{j}^{\prime}\notin E(G),$ for otherwise vertices
$w,w_{i},w_{j},w_{i}^{\prime}$ and $w_{j}^{\prime}$ induce a cycle $C_{5},$ a
contradiction too.

(iii) No vertex $x$ in
$A$ is adjacent to both $w_{i}^{\prime}$ and $w_{j}^{\prime},$ for
otherwise $w,w_{i},w_{j},w_{i}^{\prime},w_{j}^{\prime}$ and $x$ induce a cycle
$C_{6},$ a contradiction.

Accordingly, since $\delta\geq2$, each vertex $A-\{w\}$ still has a
neighbor in $A.$ In particular, $A-\{w,w_{1}^{\prime
},w_{2}^{\prime},...,w_{\Delta}^{\prime}\}$ is non empty and induce an
isolate-free subgraph. Now, consider the partition $\Phi=\{V_{1}%
,V_{2},.,V_{\Delta},V_{\Delta+1}\},$ where for any $i\in\{1,...,\Delta\},$
each $V_{i}=\{w_{i},w_{i}^{\prime}\}$ and $V_{\Delta+1}=A-\{w'_1,w'_2,...,w'_\Delta\}$. 
Clearly since $w\in V_{\Delta+1}$ and
$w$ has no neighbor in $V_{\Delta+1},$ no set of $\Phi$ is a TRD-set.
Moreover, it is not hard to notice that $V_{\Delta+1}$ forms a total
restrained coalition with any other set of $\Phi,$ leading to $C_{tr}%
(G)\geq\left\vert \Phi\right\vert =\Delta+1.$ $\Box$\newline

The bound established in Theorem 2.8 is tight, as demonstrated, for example, 
by any cycle $C_n$ where $n \not\equiv 0 \pmod{4}$ and $n \geq 7$. (see Th.~\ref{cn})

\medskip
We next present a technical lemma, which gives us the number of total restrained coalitions 
involving any set in a $C_{tr}(G)$-partition of $G$. 

\begin{lemma}\label{4}
	If $G$ is an isolate-free graph, then for any $C_{tr}(G)$-partition $\Phi$ and for any 
	$X\in \Phi$, the number of total restrained coalitions formed by $X$ is at most $\Delta(G)$.	    
 
\end{lemma} 
\begin{proof} Since $X\in\Phi$, $X$ is not a TRD-set. We now distinguish two 
cases.

\nt {\bf Case 1.} There is a vertex $v \in V(G)$ such that  $N(v) \cap X=\emptyset$.\newline
We first assume that $v\in X$. If a set $A\in \Phi$ forms a total restrained 
coalition with $X$, then $A\cup X$ is a TRD-set of $G$. So $v$ must has at 
least one neighbor in $A$. Thus, there are at most $|N(v)|-1\leq \Delta(G)-1$ 
other sets that can be in a total restrained coalition with $X$, and 
consequently, $X$ is in at most $\Delta(G)$ total restrained coalitions. 
Next let $v \not\in X$ and $X\cap N(v)=\emptyset$. Then, each set of $\Phi$ 
which is in a total restrained coalition with $X$ must contain at least 
one of the members of $N[v]$. We claim that there is no set $Y\in \Phi$ that 
forms a total restrained coalition with $X$ and $Y\cap N[v]=\{v\}$. Suppose to
the contrary that there is a set $Y\in \Phi$ that forms a total restrained 
coalition with $X$ and $Y\cap N[v]=\{v\}$. Thus $X\cup Y$ is a TRD-set. This 
implies that $v$ has a neighbor in $X\cup Y,$ contradicting our assumption 
$X\cap N(v)=\emptyset$ and $Y\cap N(v)=\emptyset$. This proves the claim. Consequently, 
there exists a unique set $Y$ among all sets of $\Phi$ forming a total restrained coalition with 
$X$, where $v$ belongs to $Y$ and $Y$ has a non-empty intersection with $N(v)$. This implies 
that the largest possible number of sets in $\Phi$ forming a total restrained coalition with $X$ is 
no more than $|N(v)|$. Therefore, the total number of sets of $\Phi$ forming a total restrained 
coalition with $X$ is at most $\Delta(G)$. 

\nt {\bf Case 2.}  There is a vertex $v \in V-X$ such that $N(v) \cap (V-X)=\emptyset$.

In this case, we prove that there is exactly one set in $\Phi$ that forms 
a total restrained coalition with $X$. Assume that $W\in\Phi\setminus \{X\}$
such that $\{X,W\}$ is a tr-coalition. If $v\not\in W$ then $v\not\in X\cup W$ and
therefore $N(v)\cap V\setminus \left(X\cup W\right)\neq \emptyset$ because
$X\cup W$ is a TRD-set in $G$. The latter is a contradiction because $N(v)\subseteq
X.$  Hence, it must be that $v\in W.$ and thus, $W$ is the only set that forms a 
total restrained coalition with $X$.

It follows from the two cases above that $X$ belongs to, at most,
 $\Delta(G)$ total restrained coalitions. $\Box$

\end{proof}

\medskip

Now we prove the following lemmas for graphs with leaves.

\begin{lemma}~\label{5}
Let $G$ be a graph with $\delta(G)=1$, and let $x$ be a leaf of $G$ and $y$ 
be the support vertex of $x$. Let $\Phi$ be a $C_{tr}(G)$-partition, and 
let $X, Y\in \Phi$ such that $x\in X$ and $y\in Y$ (possibly $X=Y$). For 
any two sets $A,B\in \Phi$ that form a total restrained coalition, we 
have $A\in \{X, Y\}$ or $B\in\{X,Y\}$. 
\end{lemma}

\begin{proof}
Since $A$ and $B$ form a total restrained coalition, $A\cup B$ is a TRD-set of $G$. 
If $A\not\in \{X,Y\}$ and $B\not\in \{X,Y\}$, then the vertex $x$ has no neighbor 
in $A\cup B$, which is a contradiction. Therefore, $A\in \{X, Y\}$ or $B\in\{X,Y\}$. 
$\Box$ \medskip
\end{proof}

\begin{remark}~\label{5b}
Since, by the definition of a total restrained dominating $S$ set, we may deduce 
that $deg(v)\ge 2$ for every vertex $v\not\in S.$ Consequently, any leaf of $G$ must
belong to $S.$
\end{remark}

We establish next an upper bound on the total restrained coalition in terms 
 of the maximum degree of $G$.

\begin{theorem}\label{6}
Let $G$ be an isolate-free graph with $\delta(G)=1$. Then, 
$C_{tr}(G) \leq \Delta(G)+1$.
\end{theorem} 

\begin{proof} 
 Let $x$ be a vertex of $G$ with $\deg(x)=1$ and let 
$\Phi=\{V_1,V_2,\ldots,V_k\}$ be a $C_{tr}(G)$-partition. Without loss of
generality, we can assume that $x\in V_1.$ If $\{V_i,V_j\}\subseteq \Phi$ form 
a total restrained coalition then, by Remark \ref{5b}, we have that 
$x\in V_i\cup V_j$. Consequently, $V_1\in \{V_i, V_j\}$. By Lemma~\ref{4}, 
$V_1$ is in total restrained coalition with at most $\Delta(G)$ sets of 
$\Phi$. Hence, $C_{tr}(G)\leq \Delta(G)+1$.
\end{proof} 

\medskip

Let us point out that the bound given by Theorem~\ref{6} is sharp. To see this, 
it is sufficient to consider the graph depicted in Figure~\ref{fig1}, where $V_1$ 
forms a tr-coalition with any of the remaining sets $V_2,V_3,$ or $V_4.$

\begin{figure}[t!]
\begin{center}
\begin{tikzpicture}[scale=0.6]

		\node [style={black_v},label=above left:{\large $v_1$}] (0) at (-7, 4) {};
		\node [style={black_v},label=above left:{\large $v_2$}] (1) at (-7, -1) {};
		\node [style={black_v},label=below left:{\large $v_3$}] (2) at (-4, 2) {};
		\node [style={black_v},label=below left:{\large $v_4$}] (3) at (-4, -3) {};
		\node [style=black_v,label=above left:{\large $v_5$}] (4) at (-1, 4) {};
		\node [style=black_v,label=above left:{\large $v_6$}] (5) at (-1, -1) {};
		\node [style={black_v},label=above right:{\large $v_7$}] (6) 
		at (2, 1.5) {};
		\node [style={black_v},label=above right:{\large $v_8$}] (7) 
		at (5.5, 1.5) {};
		\node [style=none,label=above:{\large $\Phi=\{ V_1=\{v_1,v_2\},$}] (8) 
		at (-9, -6.5) {};
		\node [style=none,label=above:{\large $V_2=\{v_3,v_4\},$}] (9) 
		at (-3, -6.5) {};
		\node [style=none,label=above:{\large $V_3=\{v_5,v_6\},$}] (10) 
		at (2, -6.5) {};
		\node [style=none,label=above:{\large $V_4=\{v_7,v_8\}\}$}] (11) 
		at (7, -6.5) {};

		\draw (0) to (4);
		\draw (4) to (2);
		\draw (2) to (0);
		\draw (0) to (1);
		\draw (1) to (3);
		\draw (3) to (2);
		\draw (3) to (5);
		\draw (5) to (1);
		\draw (4) to (6);
		\draw (6) to (5);
		\draw (6) to (7);
\end{tikzpicture}
\end{center}
\caption{A graph attaining the bound given by Theorem \ref{6}.}\label{fig1}
\end{figure}

\medskip

\begin{theorem}
\label{delta2}
Let $G$ be an isolate-free graph with $\delta(G)=2$. Then, $C_{tr}(G) 
\leq 2 \Delta(G)$.
\end{theorem}

\begin{proof}
Let $x$ be a vertex of $G$ with $\deg(x)=2$, and suppose that $N(x)=\{y,z\}$. 
Let $\Phi$ be a $C_{tr}(G)$-partition. We now distinguish the following cases.
\begin{itemize}
\item{\bf Case 1.} There is a set $U\in \Phi$ such that $\{x,y,z\}\subseteq U$. 
Then, each set of $\Phi\backslash U$ must form a total restrained coalition with $U$.
Otherwise, we would have two distinct sets $A, B\in \Phi$ forming a total restrained
coalition. Thus, $x$ must have at least one neighbor in $A \cup B$, contradicting our
supposition that $\deg(x)=2$. Therefore, by Lemma \ref{4}, $U$ is in total restrained 
coalitions with at most $\Delta(G)$ sets. Consequently, 
$C_{tr}(G)\leq \Delta(G)+1\leq 2\Delta(G)+1$.
\item{\bf Case 2.} Assume that $X, A\in \Phi$ such that $x\in X$ and 
$\{y,z\}\subseteq A$. Since $N(x)\subseteq A$, there is no set $B\neq A$ 
that forms a total restrained coalition with~$X$. So $X$ forms a total 
restrained coalition only with $A$. Moreover, $A$ does not form a total 
restrained coalition with any other set in $\Phi$ other than $X$. Otherwise, 
we would have a set $C\in \Phi$ forming a total restrained coalition with $A$. 
Thus, $x$ must have at least one neighbor outside in $A \cup C$, contradicting 
our supposition that $\deg(x)=2$. Hence, $C_{tr}(G)\leq 2$.
\item{\bf Case 3.} Assume that $Y, B\in \Phi$ such that $y\in Y$ and 
$\{x,z\}\subseteq B$. Then, each set of $\Phi\backslash\{Y,B\}$ form a 
total restrained coalition with $Y$ or $B$. Otherwise, we would have two 
distinct sets $C, D\in \Phi$ forming a total restrained coalition. Thus, 
$x$ must have at least one neighbor in $C \cup D$, contradicting our 
supposition that $\deg(x)=2$. If $Y$ and $B$ form a total restrained coalition,
by Lemma \ref{4}, we have $C_{tr}(G)\leq \Delta(G)-1+\Delta(G)-1+1+1=2\Delta(G)$.
Next, suppose that $Y$ and $B$ do not form a total restrained coalition. 
We consider two subcases.

\item {\bf Subcase 3.1.} There exists a vertex $w\in V(G)$ having no 
neighbor in $Y\cup B$. Since any set of $\Phi\backslash\{Y, B\}$ form 
a total restrained coalition with $Y$ or $B$, in order to totally restrained 
dominate the vertex $w$, any set of $\Phi\backslash\{Y, B\}$ must contain at 
least one of the members of $N(w)$. So, by Lemma \ref{4},
$C_{tr}(G)\leq |N(w)|+2\leq \Delta(G)+2\leq2\Delta(G)+1$.

\item {\bf Subcase 3.2.} There exists a vertex $w \in (V-(Y \cup B))$ such 
that $N(w) \cap (V-(Y \cup B))=\emptyset$. It follows that 
$N(w)\subseteq (Y\cup B)$. Then all TRD-sets must contain the vertex $w$, 
as each set of $\Phi\backslash\{Y,B\}$ form a total restrained coalition 
with $Y$ or $B$. This yields that $w$ is totally restrained dominated. 
Since $x$ and $y$ are adjacent, we deduce that there are at most 
$|N(y)|-1\leq \Delta(G)-1$ sets containing a member of $N(y)$. Thus, the 
set $Y$ is in at most $|N(y)|-1\leq \Delta(G)-1$ total restrained coalitions.
Analogously, we observe that the set $B$ is in at most $|N(z)|-1\leq \Delta(G)-1$ 
total restrained coalitions. Hence, 
$C_{tr}(G)\leq \Delta(G)-1 + \Delta(G)-1+2=2\Delta(G) \leq 2\Delta(G)+1$.

\item {\bf Case 4.} There are two distinct sets $Z, C\in \Phi$ such 
that $z\in Z$ and $\{x,y\}\subseteq C$.
The proof is similar to the proof of  {\bf Case 3}.

\item {\bf Case 5.} 
Assume that $X, Y, Z\in \Phi$ such that $x\in X, y\in Y$ and $z\in Z$.
We claim the following facts,
\begin{itemize}
	\item[(5.i)] If $X,T \in \Phi$ form a tr-coalition then $T\in\{Y,Z\}$. This is because
	the neighbors of $x$ belongs to $Y\cup Z.$
	\item[(5.ii)] $Y,Z$ can not form a tr-coalition because otherwise $x\not\in Y\cup Z$ would not be
	total restrained dominated.
	\item[(5.iii)] If $Y,T\in \Phi\setminus \{X,Z\}$ form a tr-coalition then 
	$N(z)\cap \left(Y\cup T\right)\neq\emptyset.$
	Otherwise,
	the vertex $z$, which does not belongs to $Y\cup T,$ would not be total restrained dominated by $Y\cup T$. 
	\item[(5.iv)] If $Z,T\in \Phi\setminus \{X,Y\}$ form a tr-coalition then 
	$N(y)\cap \left(Z\cup T\right)\neq\emptyset.$
	Otherwise,
	the vertex $y$, which does not belongs to $Z\cup T,$ would not be total restrained dominated by $Z\cup T$. 
\end{itemize}

Now, let us distinguish three different cases,

\begin{itemize}
	\item If $N(z)\cap Z \neq \emptyset$ or $N(z) \cap Y\neq \emptyset$ then by considering (5.iii) we know that
	$Y$ can form a tr-coalition with, at most, $|N(z)|-2$ different sets $T$. Since $x$ and $y$ are adjacent, we 
	deduce that there are at most $|N(y)|-1\leq \Delta(G)-1$ sets which contain a member of $N(y)$. 
	Thus, the set $Z$ is in at most $|N(y)|-1\leq \Delta(G)-1$ total restrained coalitions. Therefore, 
	$$ C_{tr}(G)\le |N(z)|-2+|N(y)|-1+3\le 2\Delta(G)$$
	\item If $N(y)\cap Z \neq \emptyset$ or $N(y) \cap Y\neq \emptyset$ then by considering (5.iv) we know that
	$Z$ can form a tr-coalition with, at most, $|N(y)|-2$ different sets $T$. Since $x$ and $z$ are adjacent, we 
	deduce that there are at most $|N(z)|-1\leq \Delta(G)-1$ sets which contain a member of $N(z)$. 
	Thus, the set $Y$ is in at most $|N(z)|-1\leq \Delta(G)-1$ total restrained coalitions. Therefore, 
	$$ C_{tr}(G)\le |N(z)|-1+|N(y)|-2+3\le 2\Delta(G)$$
	\item Otherwise, assume that $N(z)\cap Z=N(z)\cap Y=N(y)\cap Z=N(y)\cap Y = \emptyset.$ 
	If $T$ form a tr-coalition with $Y$ then $N(z)\cap T\neq\emptyset$ because $z\not\in Y\cup T$ and
	$Y\cup T$ is a TRD-set. Besides, $N(y)\cap T\neq\emptyset$ because $y\in Y\cup T$, $N(y)\cap Y
	=\emptyset$ and $Y\cup T$ is a TRD-set. Consequently, any set $T$ that forms a tr-coalition with $Y$ 
	(analogously, with $Z$) must contain both a neighbor of $y$ and a neighbor of $z$. Therefore,
	$$ C_{tr}(G)\le |N(z)|-1+3\le \Delta(G)+2\le 2\Delta(G).$$
\end{itemize}

Based on the analysis of all the above cases, we infer that $C_{tr}(G)\le 2\Delta(G).$ 
$\Box$

\end{itemize}
\end{proof}

The bound described in Theorem~\ref{delta2} is sharp, as illustrated by
any cycle \( C_n \) with \( n \geq 7 \) and \( n \equiv 0 \pmod{4} \)
(refer to Th.~\ref{cn} for further details).


\section{ Total restrained coalition number of specific graphs }

In this section, we deal with the problem of obtaining the exact value 
of the  total restrained coalition number. We first recall the following results.

\begin{proposition}{\rm\cite{A6}}
Let $n \geq 4$ be a positive integer. Then $\gamma_{tr}(K_n)=2$.
\end{proposition}

\begin{proposition} {\rm\cite{A6}}\label{7} 
Let $n_1$ and $n_2$ be positive integers such that $\min\{n_1, n_2\} \geq 2$. 
Then $\gamma_{tr}(K_{n_1, n_2})=2$.
\end{proposition}

\begin{proposition} {\rm\cite{A6}} \label{8}
Let $n$ be a positive integer. Then $\gamma_{tr}(K_{1,{n-1}})=n$.
\end{proposition}

The following proposition gives us the total restrained coalition number of 
the complete graph.

\begin{proposition} \label{9}
	Let $n \geq 4$ be a positive integer. Then $C_{tr}(K_n)=n$.
\end{proposition} 

\begin{proof} Let $G$ be a complete graph of order $n$ with vertex set 
$V=\{v_1, v_2,\ldots, v_n\}$. Since $\gamma_{tr}(G)=2$, every two adjacent 
vertices $v_i$ and $v_j$ of $G$ can be in a total restrained coalition. 
It follows that $\Phi=\left\{\{v_1\}, \{v_2\}, \ldots, \{v_n\}\right\}$ 
is a trc-partition, and hence $C_{tr}(K_n)=n$. $\Box$ 

\end{proof}\medskip

By Proposition \ref{7}, we get the following result.
\begin{observation} \label{10}
Let $G=K_{p,q}$ be a complete bipartite graph such that $q\geq p\geq 2$. 
Then $C_{tr}(K_{p,q})=p+q=n$.

\end{observation}

Proposition \ref{8} gives the next result.

\begin{observation} \label{11}
If $G=K_{1,{n-1}}$ is a star graph, then $C_{tr}(K_{1,{n-1}})=2$.

\end{observation}

Next we determine the total restrained coalition number of paths. But before we 
need to recall the following result from \cite{A6}.
 
 \begin{theorem}{\rm\cite{A6}} \label{12}
 Let $n\geq 4$ be a positive integer. Then 
 $\gamma_{tr}(P_n)=n-2\lfloor\frac{n-2}{4}\rfloor$.
 \end{theorem}

\begin{theorem}\label{13}
For any path $P_n$,
		\begin{equation*}
		C_{tr}(P_n)=\left\{
		\begin{aligned}[c]
		2, & { \ \ if\ } 2\leq n\leq 7 \\[.5em]
		3, & { \ \ if \ } n\geq 8 \\
		\end{aligned}\right.
		\end{equation*}
	\end{theorem}

\begin{proof} Let $V(P_n)=\{v_1,v_2,\ldots,v_n\}$. By Theorem \ref{6} and 
Corollary \ref{3}, we have $2\leq C_{tr}(P_n)\leq 3$ for any path $P_n$. 
If $n=2$, then Proposition \ref{9} gives the desired result, while if $n=3$, 
the result follows from Observation \ref{11}. We next proceed to show that 
$C_{tr}(P_n) \ne 3$ where $4 \leq n\leq 7$. Let $\Phi=\{A, B, C\}$ be a 
$C_{tr}(P_n)$-partition. By Lemma \ref{4}, each set of $\Phi$ is in total 
restrained coalition with at most two sets of $\Phi$. So, without loss of 
generality, assume that each of $B$ and $C$ forms a total restrained coalition 
with $A$. By Theorem \ref{12}, we have $|A|+|B|\geq n-2\lfloor\frac{n-2}{4}\rfloor$ 
and $|A|+|C|\geq n-2\lfloor\frac{n-2}{4}\rfloor$. Therefore, 
$2|A|+|B|+|C|\geq 2n-4\lfloor\frac{n-2}{4}\rfloor$. On the other hand, we know
that $|A|+|B|+|C|=n$. Hence, $|A|\geq n-4\lfloor\frac{n-2}{4}\rfloor$. 
Now suppose that $n=4$. Hence, $|A|\geq 4$, contradicting the fact that 
$|A|<4$. This implies that $C_{tr}(P_4)\neq 3$. If $n=5$, then $|A|\geq 5$ 
which impossible as $|A|<5$. Consequently, $C_{tr}(P_5)\neq 3$. Now assume 
that $n=6$. Thus, we have $|A|\geq 2$. On the other side, $|A|\leq 5$. We 
now distinguish the following cases.

   \nt {\bf Case 1.} $\Phi$ consists of a set of cardinality 2 (namely $A$), 
a set of cardinality 3 (namely $B$) and a singleton set (namely $C$). 
Since $\gamma_{tr}(P_6)=4$, each of $A$ and $C$ must be in a total 
restrained coalition with $B$. This is impossible because $P_6$ has no 
TRD-set of order 5. Hence, $C_{tr}(P_6) \ne 3$.

   \nt {\bf Case 2.} Let $|A|=|B|=|C|=2$. We may assume that each of $B$ and 
$C$ must be in a total restrained coalition with $A$, which is impossible, 
as $P_6$ has a unique TRD-set of order 4. Hence, $C_{tr}(P_6) \ne 3$.

   \nt {\bf Case 3.} $\Phi$ consists of a set of cardinality 3 (namely $A$), a 
   set of cardinality 2 (namely $B$) and a singleton set (namely $C$). 
   Analogous argument as in Case 1(by interchanging the roles of $A$ and $B$) 
   can be applied to show that $\Phi$ of order 3 does not exist.

\nt {\bf Case 4.} $\Phi$ consists of a set of cardinality 4, say $A$, and two 
singleton sets such as $B$ and $C$. Since $\gamma_{tr}(P_6)=4$, no two singleton 
sets in $\Phi$ form a total restrained coalition. It follows that each of $B$ 
and $C$ must be in a total restrained coalition with $A$, which is impossible, 
as $P_6$ has no TRD-set of order 5. Hence, $C_{tr}(P_6) \ne 3$.

\nt {\bf Case 5.} Let $|A|=5$. It follows that either $B$ or $C$ is an empty set. 
But this is impossible. Then, $C_{tr}(P_6) \ne 3$.

Next suppose that $n=7$. So, we have $|A|\geq 3$. On the other side, $|A|\leq 6$. 
We now consider the following cases.

\nt {\bf Case 1.} $\Phi$ consists of two sets of cardinality 3, say $A$ and $B$, 
and a singleton set $C$. Since $\gamma_{tr}(P_7)=5$, neither $A$ nor $B$ can be 
in a total restrained coalition with $C$. Consequently, there is no total 
restrained coalition partition of order 3. Hence, $C_{tr}(P_7) \ne 3$.

\nt {\bf Case 2.} $\Phi$ consists of a set of cardinality 4 (namely $A$), a 
set of cardinality of 2 (namely $B$) and a singleton set (namely $C$). Since 
$\gamma_{tr}(P_7)=5$, each of $B$ and $C$ must be in a total restrained 
coalition with $A$. This is impossible because $P_7$ has no TRD-set of order 6. 
Thus, $C_{tr}(P_7) \ne 3$.

\nt {\bf Case 3.} $\Phi$ consists of a set of cardinality 5, say $A$, and 
two singleton sets such as $B$ and $C$. Since $\gamma_{tr}(P_7)=5$, each 
of $B$ and $C$ must be in a total restrained coalition with $A$. This is 
impossible because $P_7$ has no TRD-set of order 6. Hence, $C_{tr}(P_7) \ne 3$.

\nt {\bf Case 4.} Let $|A|=6$. It follows that either $B$ or $C$ is an empty set. But this
is impossible. Hence, $C_{tr}(P_7) \ne 3$.\medskip

By the above discussions, we infer that $C_{tr}(P_n)=2$ where $4 \leq n\leq 7$.\medskip

Finally, let $n\geq 8$. By Theorem \ref{6}, for any path $P_n$ we have 
$C_{tr}(P_{n})\leq 3$. To achieve equality, all we need is to give a total 
restrained partition of order 3 for any $n\geq8$, and which will be as follows:  

$$\Phi(P_n)= \left\{X=\{v_1,v_2 \dots v_{n-6},v_{n-1},v_{n}\}, 
Y=\{v_{n-5},v_{n-4}\}, Z=\{v_{n-3},v_{n-2}\}\right\}.$$
  
  One can observe that each of $Y$ and $Z$ is in a total restrained coalition  
  with $X$. Therefore, the proof is complete. $\Box$

\end{proof}\medskip

We close this section by calculating the total restrained
coalition number of cycles. It is straightforward to see that $C_{tr}(C_3)=2$
so we next focus on the cases where the order is at least $4$.
We begin by recalling the following result.

 \begin{theorem}{\rm\cite{A6}} \label{cn}
 Let $n\geq 4$ be a positive integer. Then 
 $\gamma_{tr}(C_n)=n-2\lfloor\frac{n}{4}\rfloor$.
 \end{theorem}

\begin{theorem}\label{cycle}
	Let $n\ge 4$ be an intenger and let $C_n$ be a cycle. Then, 
	$$ \mathcal{C}_{tr}(C_{n})=\left\{
	\begin{array}{cc}
	4  &\quad n\equiv 0 ~(\mbox{mod } 4)\\[.5em]
	3   &\quad\mbox{otherwise. }
	\end{array}\right.
	$$
\end{theorem}  

\begin{proof}  Let $G=C_n$, where $V(G)=\lbrace v_1,v_2,\dots ,v_n\rbrace$ and 
$E(G)=\lbrace v_i v_{i+1} : 1\leq i\leq n\rbrace$. By Theorem \ref{delta2} and 
Corollary \ref{3}, we have $2\leq C_{tr}(C_n)\leq 4$ for any cycle $C_n$. 
Now assume that $n\equiv 0~(\mbox{mod } 4)$. 
We find a total restrained coalition partition of order 4 as follows:

$\Phi(C_n)=\{A=\{v_{n}\}, B=\{v_{n-2}\}, C=\{v_1,v_2,v_5,v_6 \dots v_{n-7},v_{n-6}, v_{n-3}\}, 
D=\{v_3,v_4,v_7,v_8 \dots v_{n-5},v_{n-4}, v_{n-1}\} \}$. It is easy to verify that $A$ and 
$D$ form a total restrained coalition, and $B$ and $C$ form a total restrained coalition.

We next suppose that $n$ is not divisible by 4. We claim that 
$\mathcal{C}_{tr}(C_n)\neq 4$. Suppose, to the contrary that 
$\mathcal{C}_{tr}(C_n)=4$. Let $\{A, B, C, D\}$ be a 
$\mathcal{C}_{tr}(C_n)$-partition.  
	By Theorem~\ref{4}, we conclude that each set of $\Phi$ forms a total 
	coalition with at most two sets of $\Phi$. Now assume without loss of 
	generality that $A$ and $B$ form a total restrained coalition, and $C$ 
	and $D$ form a restrained total coalition. By Theorem~\ref{cn}, we have 
	$\gamma_{tr}(C_n)=n-2\lfloor\frac{n}{4}\rfloor$ where $n\geq 4$. Since 
	$A \cup B$ and $C \cup D$ are the total restrained dominating sets, it 
	follows that $|A|+|B|\geq n-2\lfloor\frac{n}{4}\rfloor$ and $|C|+|D|\geq
	 n-2\lfloor\frac{n}{4}\rfloor$. Hence, $|A|+|B|+|C|+|D|\geq 2n-4\lfloor
	 \frac{n}{4}\rfloor$. On the other side, we know that $|A|+|B|+|C|+|D|=n$. 
	 Thus, we obtain $4\lfloor\frac{n}{4}\rfloor \geq n$. Since $n$ is not 
	 divisible by 4, we get $n-1\geq n$. But this leads to a contradiction. 
	 Then, $\mathcal{C}_{tr}(C_n)\neq 4$. Hence, $\mathcal{C}_{tr}(C_n)\leq 3$. 
	 In the following, we construct a total restrained coalition partition 
	 of order 3. So assume that $\Phi =\{A=\{v_1, v_2 \dots v_{n-4}\}, 
	 B=\{v_{n},v_{n-1}\}, C=\{v_{n-2},v_{n-3}\}\}$. Note that each of $B$ 
	 and $C$ forms a total restrained coalition with $A$. $\Box$
    
\end{proof}

\section{Graphs with large total restrained coalition number}
In this section, we will be interested in connected graphs $G$ of order $n$ such that 
$C_{tr}(G)=n$. We start by giving the following sufficient condition for $G$ to have 
$C_{tr}(G)=n$. We recall that a universal vertex of a graph $G$ is a vertex that is 
adjacent to every other vertex in $G$.

\begin{proposition}
If $G$ is a graph of order $n$ and minimum degree at least three and having a
universal vertex, then $C_{tr}(G)=n$.
\end{proposition}

\textbf{Proof. }Let $v_{1},v_{2},...,v_{n}$ denote the vertices of $G,$ and
assume that $v_{1}$ is a universal vertex in $G.$ Clearly, $n\geq4$, since $\delta(G)\geq3.$ 
Consider the partition $\Phi=\{\{v_{1}\},\{v_{2}%
\},...,\{v_{n}\}\}.$ Since each vertex in $G$ has at least two neighbors
besides $v_{1},$ the set $\{v_{1}\}$ forms a total restrained coalition with any other
set of $\Phi.$ Hence $\Phi$ is a trc-partition in $G,$ and thus $C_{tr}%
(G)=\left\vert \Phi\right\vert =n.$ $\Box$ \medskip

We note that the condition on the minimum degree to be at least three is important, and this 
can be seen for stars of order at least three, when $\delta(G)=1$, and for the graph $G$ 
obtained from $k\geq2$ triangles sharing the same vertex, when $\delta(G)=2$. Both graphs have 
a universal vertex but in either graph $C_{tr}(G)<n$.

\bigskip

We now give a necessary condition for connected graphs $G$ such that
$C_{tr}(G)=n.$

\begin{proposition}
\label{Ctr=n}Let $G$ be a connected graph of order $n\geq3$ such that
$C_{tr}(G)=n$. Then
\begin{enumerate}
\item[(1)] $\gamma_{tr}(G)=2.$

\item[(2)] $\delta(G)\geq2.$

\item[(3)] $G$ has diameter at most $2.$
\end{enumerate}
\end{proposition}

\textbf{Proof. }Let\textbf{ }$\Phi=\{\{v_{1}\},\{v_{2}\},...,\{v_{n}\}\}$ be a
$C_{tr}(G)$-partition. Clearly, item (1) follows from the fact that each set
of\ $\Phi$ has to be in a total restrained coalition with another set of
$\Phi.$

To prove item (2), suppose to the contrary that $G$ contains a leaf, say
$v_{1},$ and let $v_{2}$ be the neighbor of $v_{1}.$ Since $v_{1}$ and $v_{2}$
belong to all TRD-sets of $G$, and since $n\geq3,$ it is no longer possible to
form total restrained coalitions from any two singleton sets of $\Phi,$ contradicting
$\gamma_{tr}(G)=2.$ Hence $\delta(G)\geq2.$ Finally, to prove item (3), we
first note that by (1) $G$ has diameter at most three. Assume, for a
contradiction, that $G$ has diameter three, and let $x$ and $y$ be two
vertices at distance three in $G.$ Observe that $x$ and $y$ have no common
neighbor, and thus set $\{x\}$ cannot form a total restrained coalition with
any other set of $\Phi$ so that $\gamma_{tr}(G)=2.$ Hence $G$ has diameter at
most $2.$ $\Box$

It is worth noting that the converse of Proposition \ref{Ctr=n} is not true. This can be 
seen by the graph $G$ obtained from a cycle $C_5$ whose vertices are labeled in order 
$v_1,v_2,v_3,v_4,v_5,v_1$ by adding the edge $v_2v_5$. Clearly $G$ satisfies conditions 
of Proposition \ref{Ctr=n}, but $\{v_1\}$ cannot form a total restrained coalition with 
any singleton set of $V$, leading to $C_{tr}(G)<5$.

\bigskip

In the next, we characterize connected triangle-free graphs $G$ of order
$n\geq2$ with $C_{tr}(G)=n$.

\begin{proposition}
\label{triangle-free}A connected triangle-free graphs $G$ of order $n\geq2$
satisfies $C_{tr}(G)=n$ if and only if $G=P_{2}$ or $G=K_{p,q}$ for any
integers $p,q\geq2.$ 
\end{proposition}

\textbf{Proof. }Let $G$ be a connected triangle-free graphs of order $n\geq2$
such that $C_{tr}(G)=n.$ Clearly, if $n=2,$ then $G$ is a path $P_{2}.$ Hence
assume that $n\geq3.$ By item (1) of Proposition \ref{Ctr=n}, $\gamma
_{tr}(G)=2,$ and so let $\{u,v\}$ be a minimum TRD-set of $G.$ Since $G$ is
triangle free, each of $N(v)$ and $N(v)$ is an independent set. Moreover, $u$
and $v$ have no common neighbor. Let $A=N(u)-\{v\}$ and $B=N(v)-\{u\}.$
Observe that if $A=\emptyset,$ then $G$ is a star and thus $C_{tr}(G)=2<n,$ 
{  a contradiction}. Therefore $A\neq\emptyset$ and likewise 
$B\neq\emptyset.$ Now,
since $\delta(G)\geq2,$ by Proposition \ref{Ctr=n}-(2), each vertex in $A$ has
at least one neighbor in $B$ {  and vice versa}.
 Moreover, if a vertex $x\in A$
has a non-neighbor $y\in B,$ then clearly $x$ and $y$ are at distance three in
$G,$ leading that $G$ has diameter 3, contradicting item (3) of Proposition
\ref{Ctr=n}. Hence $A\cup B$ induces a complete bipartite graph $K_{\left\vert
A\right\vert ,\left\vert B\right\vert },$ and therefore $G=K_{\left\vert
A\right\vert +1,\left\vert B\right\vert +1},$ as desired. 

For the converse, since any two adjacent vertices of either $P_{2}$ or
$K_{p,q}$ with $p,q\geq2,$ form a TRD-set of $G,$ we easily conclude that
$C_{tr}(G)=n.$ $\Box$

\bigskip

According to Propositions \ref{Ctr=n} and \ref{triangle-free}, for every tree
$T$ of order $n\geq3,$ $C_{tr}(T)\leq n-1.$ Our aim in the following is to
characterize all trees $T$ of order $n\geq3$ such that $C_{tr}(T)=n-1.$

\begin{theorem}
The path $P_{3}$ is the only tree $T$ of order $n\geq3$ satisfying
$C_{tr}(T)=n-1.$ 
\end{theorem}

\textbf{Proof. }If $T=P_{3},$ then by Theorem \ref{13}, $C_{tr}(P_{3})=2=n-1.$
To prove the converse, let $T$ be a tree of order $n\geq3$ such that
$C_{tr}(T)=n-1.$ Since $C_{tr}(T)\leq\Delta(T)+1,$ by Theorem \ref{6}, we
deduce that $\Delta(T)\geq n-2.$ If $\Delta(T)=n-1,$ then $T$ is a star
$K_{1,n-1},$ and thus by Observation \ref{11}, $n=3$ leading to $T=K_{1,2}%
=P_{3}.$ Hence in the following we can assume that $\Delta(T)=n-2,$ and thus
$n\geq4.$ Clearly, $T$ is a double star $S_{n-3,1},$ with one support vertex
having $n-3$ leaf neighbors and the other with only one leaf neighbor.
Let\textbf{\ }$\Phi$ be a $C_{tr}(T)$-partition.\ Since $\left\vert
\Phi\right\vert =n-1,$ one member of $\Phi$ has cardinality 2, and any other
member is a singleton set, that is of cardinality 1. Now, if $n\geq5,$ then
some leaf $w$ of $T$ alone will form a singleton set $\{w\}$ in $\Phi.$ But
then $\{w\}$ cannot be in total restrained coalition with any other set of
$\Phi,$ a contradiction. Hence $n=4$ and so $T$ is simply a path $P_{4}.$ By
Theorem \ref{13}, $C_{tr}(P_{4})=2<n-1,$ a contradiction. This completes the
proof. $\Box$

\bigskip

From the above, the following corollary is derived.

\begin{corollary}
If $T$ is a tree of order $n\geq4,$ then $C_{tr}(T)\leq n-2.$
\end{corollary}

\begin{proposition}
Let $T$ be a tree of order $n\geq4.$ Then $C_{tr}(T)=n-2$ if and only if
$T\in\{K_{1,3},P_{4},P_{5}\}$
\end{proposition}

\textbf{Proof. }Let $T$ be a tree of order $n\geq4$ such that $C_{tr}(T)=n-2.$
By Theorem \ref{6}, $C_{tr}(T)\leq\Delta(T)+1$ and thus $\Delta(T)\geq n-3.$
If $\Delta(T)=n-1,$ then $T$ is a star $K_{1,n-1},$ and by Observation
\ref{11} we deduce that $n=4,$ that is $T=K_{1,3}.$ Hence in the following we
only consider $\Delta(T)\in\{n-3,n-2\}.$ 

Assume first that $\Delta(T)=n-2.$ Then $T$ is a double star $S_{n-3,1},$ and
since $V(T)$ is the unique TRD-set, we deduce that any $C_{tr}(T)$-partition
$\Phi$ contains only two sets, leading to $\left\vert \Phi\right\vert =n-2=2.$
Therefore $n=4$ and $T=S_{1,1}$, that is $T$ is a path $P_{4}.$ 

Finally, assume that $\Delta(T)=n-3.$ Let $x$ be the vertex of maximum degree,
and let $u$ and $v$ denote the two non-neighbors of $x.$ If $u$ and $v$ are
adjacent, then, without loss of generality, let $u$ and $x$ have a common
neighbor $y.$ As before, since $V(T)$ is the unique TRD-set of $T,$ we deduce
that any $C_{tr}(T)$-partition $\Phi$ contains only two sets, leading to
$\left\vert \Phi\right\vert =n-3=2.$ Therefore, $n=5,$ and $T$ is a path
$P_{5}.$ In the following, we can assume that $uv\notin E(T).$ In this case,
let $y$ and $z$ be the common neighbors of $x$ with $u$ and $v,$ respectively.
Note that every neighbor of $x$ besides $y$ and $z,$ if any, is a leaf. But
since $T$ has diameter $4,$ any $C_{tr}(T)$-partition $\Phi$ contains only two
sets, leading as before to $n=5,$ and thus $T=P_{5}.$ 

The converse follows from Observation \ref{11} and Theorem \ref{13}.  $\Box$

\section{Concluding remarks}
In this article, we introduced and studied the total restrained coalition in graphs. 
We investigated the existence of a total restrained coalition partition. We derived 
some upper and lower bounds for the total restrained coalition number. We provided a
 necessary condition for graphs $G$ of order $n$ such that $C_{tr}(G)=n$, and we characterized 
 triangle-free graphs with $C_{tr}(G)=n$. Restricted to the class of trees,  we characterized 
 those trees $T$ of order $n$ such that $C_{tr}(T)$ belongs to $\{n-1,n-2\}$. In order to expand 
 the study of coalitions, we propose some potential research directions.

\begin{enumerate}
	\item Characterize all graphs satisfying $C_{tr}(G)=C(G)$.

        \item Study the total restrained coalitions in cubic graphs.
        
        \item Similar to the coalition graph, it is natural to define and study the total restrained 
        coalition graph for a given graph $G$ with respect to the total restrained coalition partition 
        $\Phi$. We define it as follows. Corresponding to any total restrained coalition partition 
        $\Phi=\{V_1,V_2,\ldots, V_k\}$ of a given graph $G$, a {\em total resinated coalition graph} 
        $TRCG(G,\Phi)$ is associated with a one-to-one correspondence between the vertices of $TRCG(G,\Phi)$ 
        and the sets $V_1,V_2,\ldots,V_k$ of $\Phi$. Two vertices of $TRCG(G,\Phi)$ are adjacent if 
        
        and only if their corresponding sets in $\Phi$ form a total restrained coalition. 

       \item  Study the restrained coalitions in graphs.

    \item Does there exist a polynomial algorithm for computing $C_{tr}(T)$ for any tree $T$?

  \item Study the complexity issue of the decision problem related to the total restrained coalition number. 

\end{enumerate}

\section*{Author contributions statement}
All authors contributed equally to this work.

\section*{Conflicts of interest}
The authors declare no conflict of interest.

\section*{Data availability}
No data was used in this investigation.

\end{document}